%% file: manuscript-numapde-preprint.tex
\pdfoutput=1

\makeatletter
\IfFileExists{./numapde-latex/numapde-packages.sty}{\providecommand*{\input@path}{}\edef\input@path{{./numapde-latex/}\input@path}}{}
\makeatother

\RequirePackage{algpseudocode}
\documentclass{numapde-preprint}

\usepackage{numapde-semantic}
\usepackage{numapde-style}
\usepackage{numapde-local}

\IfFileExists{./numapde-bibliography/numapde.bib}{\addbibresource{./numapde-bibliography/numapde.bib}}{\addbibresource{numapde.bib}}

\hypersetup{
	pdftitle={The Proximal Map of the Weighted Mean Absolute Error},
	pdfauthor={Lukas Baumgärtner, Roland Herzog, Stephan Schmidt, Manuel Weiß},
	pdfkeywords={proximal map, weighted mean absolute error, multi-thresholding}
}

\title{The Proximal Map of the Weighted Mean Absolute Error\thanks{This work was supported by DFG grants HE~6077/10--2 and SCHM~3248/2--2 within the Priority Program SPP~1962 (Non-smooth and Complementarity-based Distributed Parameter Systems: Simulation and Hierarchical Optimization), which is gratefully acknowledged.}}
\subtitle{}
\shorttitle{The Proximal Map of the WMAE}

\author{Lukas Baumgärtner\thanks{Institut für Mathematik, Humboldt-Universität zu Berlin, 10099 Berlin, Germany (\email{lukas.baumgaertner@hu-berlin.de}, \url{https://www.mathematik.hu-berlin.de/en/people/mem-vz/1693318}, \orcid{0000-0003-1007-4815}, \email{s.schmidt@hu-berlin.de}, \url{https://www.mathematik.hu-berlin.de/en/people/mem-vz/1693090}, \orcid{0000-0002-4888-0794}).}
\and
Roland Herzog\thanks{Interdisciplinary Center for Scientific Computing, Heidelberg University, 69120 Heidelberg, Germany (\email{roland.herzog@iwr.uni-heidelberg.de}, \url{https://scoop.iwr.uni-heidelberg.de}, \orcid{0000-0003-2164-6575}, \email{manuel.weiss@iwr.uni-heidelberg.de}, \url{https://scoop.iwr.uni-heidelberg.de}, \orcid{0000-0002-6098-9725}).}
\and
Stephan Schmidt\footnotemark[2]
\and
Manuel Weiß\footnotemark[3]}
\shortauthor{L. Baumgärtner, R. Herzog, S. Schmidt and M. Weiß}

\dedication{}

\IfFileExists{numapde-uncolorizeHyperrefIfChanges.sty}{\RequirePackage{numapde-uncolorizeHyperrefIfChanges}}{}

\begin{document}
\maketitle

\begin{abstract}
We investigate the proximal map for the weighted mean absolute error function. 
An algorithm for its efficient and vectorized evaluation is presented.
As a demonstration, this algorithm is applied as part of a checkerboard algorithm to solve a total-variation image denoising (ROF) problem as well as a non-smooth energy minimization problem.\end{abstract}


\begin{AMS}
\href{https://mathscinet.ams.org/msc/msc2010.html?t=90C25}{90C25}, \href{https://mathscinet.ams.org/msc/msc2010.html?t=68U10}{68U10}, \href{https://mathscinet.ams.org/msc/msc2010.html?t=94A08}{94A08}, \href{https://mathscinet.ams.org/msc/msc2010.html?t=46N10}{46N10}
\end{AMS}

\input{main.tex}

\printbibliography

\end{document}

%% file: main.tex
\section{Introduction}
\label{section:Introduction}

The proximity operator, or proximal map, plays a fundamental role in non-smooth optimization; see for instance \cite{ChambollePock:2011:1,CombettesPesquet:2011:1,ParikhBoyd:2014:1}.
Given a function $f \colon \R^n \to \R \cup \{\infty\}$, the proximal map $\prox{f} \colon \R^n \to \R^n$ is defined as the solution of the problem
\begin{equation}
	\label{eq:prox-problem}
	\text{Minimize}
	\quad
	f(y) + \frac{1}{2} \norm{y - x}_2^2
	,
	\quad
	\text{where }
	y \in \R^n
	.
\end{equation}
Under the mild condition that $f$ is proper, lower semicontinuous and convex, $\prox{f}$ is well-defined.
We refer the reader to \cite[Ch.~12.4]{BauschkeCombettes:2011:1} for details and further properties.

In this paper we present theory and an efficient algorithm for the evaluation of $\prox f$, where $f \colon \R \to \R$ is defined as
\begin{equation}
	\label{eq:weighted_mean_absolute_error}
	f(x)
	\coloneqq
	\sum_{i=1}^N \weight_i \, \abs{x - \data_i}
	.
\end{equation}
Here $\weight_i > 0$ are given, positive weights and $\data_i \in \R$ are given data, $i = 1, \ldots, N$ for some $N \in \N$.
We refer to \eqref{eq:weighted_mean_absolute_error} as the weighted mean absolute error.
Any of its minimizers is known as a weighted median of the data $\{\data_i\}$.
Clearly, $f$ is proper, continuous and convex, and so is $\gamma \, f$ for any $\gamma > 0$.

By definition, the proximal map $\prox{\gamma f} \colon \R \to \R$ for $f$ as in \eqref{eq:weighted_mean_absolute_error} is given by
\begin{equation}
	\label{eq:prox_WMAE}
	\prox{\gamma f}(x)
	\coloneqq
	\argmin_{y \in \R} \gamma \sum_{i=1}^N \weight_i \, \abs{y - \data_i} + \frac{1}{2} (y - x)^2
	.
\end{equation}
In the case of a single data point ($N = 1$), problem \eqref{eq:prox_WMAE} reduces to the well-known problem
\begin{equation}
	\label{eq:prox_WMAE_1D}
	\argmin_{y \in \R} \gamma \weight \, \abs{y - \data} + \frac{1}{2} (y - x)^2
\end{equation}
with $w > 0$ and $d \in \R$, whose unique solution is explicitly given in terms of the soft-thresholding operator $S_r(x) \coloneqq \max\{0, \, \abs{x} - r\} \sgn(x)$.
In this case, we have
\begin{equation}
	\label{eq:prox_WMAE_1D_solution}
	\prox{\gamma f}(x)
	=
	\data + S_{\gamma \weight}(x - \data)
	=
	\data + \max\paren[auto]\{\}{0, \, \abs{x - \data} - \gamma \weight} \sgn(x-\data)
	.
\end{equation}
This map, often with $d = 0$, arises in many iterative schemes for the solution of problems involving the $1$-norm; see for instance \cite{DaubechiesDefriseDeMol:2004:1,GoldsteinOsher:2009:1}.
We can therefore view \eqref{eq:prox_WMAE} as a multi-thresholding operation.

We wish to point out that our problem of interest \eqref{eq:prox_WMAE} is different from the LASSO problem
\begin{equation}
	\label{eq:LASSO}
	\text{Minimize}
	\quad
	\norm{Ay - d}_2^2 + \lambda \, \norm{y}_1
	\quad
	\text{where }
	y \in \R^n 
	,
\end{equation}
see \cite{Tibshirani:1996:1,ChenDonohoSaunders:1998:1}.
In the latter, $y$ is multi-dimensional and the deviation of its image under a linear map $A$ from a data vector $d$ is measured.
By contrast, in \eqref{eq:prox_WMAE} we measure the deviation of a scalar $y$ from multiple data points $\data_i$.
Moreover, the roles of the $1$-norm and the $2$-norm are reversed in \eqref{eq:prox_WMAE} and \eqref{eq:LASSO}.

We point out that \cite{LiOsher:2009:1} have considered the slightly more general problem 
\begin{equation}
	\label{eq:LiOsher_problem}
  \argmin_{y \in \R} \sum_{i=1}^N \weight_i \, \abs{y - \data_i} + F(y)
\end{equation}
with $F$ strictly convex, differentiable and $F^\prime$ bijective.
The prototypical examples are functions $F(y) = \lambda \, \abs{y - x}^\alpha$ with $\alpha > 1$.
We concentrate on the case $F(y) = \frac{1}{2 \gamma} (y-x)^2$, which agrees with \eqref{eq:prox_WMAE}.
In contrast to \cite{LiOsher:2009:1}, we provide a vectorized, open-source implementation of \eqref{eq:prox_WMAE}; see \cite{BaumgaertnerHerzogSchmidtWeiss:2022:2}.
We also demonstrate the utility of our implementation of \eqref{eq:prox_WMAE} by solving, similarly as in \cite{LiOsher:2009:1}, an image denoising problem using a block coordinate descent (checkerboard) algorithm.
In order to overcome the generic failure of convergence of such a method to the unique minimizer, we combine it with restarts based on the steepest descent direction.
The emphasis in this paper, however, is on the efficient solution of \eqref{eq:prox_WMAE}.

This paper is structured as follows:
We establish an algorithm for the evaluation of the proximal map of the weighted mean absolute error \eqref{eq:prox_WMAE} in \cref{section:Evaluation} and prove its correctness in \cref{theorem:proof_of_alg}. 
In \cref{section:Structure} we briefly discuss the structural properties of the proximal map.
We conclude by showing two applications of the proposed algorithm.
The first application is to an image denoising problem using the ROF model \cite{RudinOsherFatemi:1992:1}, see \cref{section:ImageDenoising}.
In \cref{section:Application_to_Deflection_Energy_Minimization}, we address a non-smooth energy minimization problem.

\section{Algorithm for the Evaluation of the Proximal Map}
\label{section:Evaluation}

In this section, we derive an efficient algorithm for the evaluation of the proximal map $\prox{\gamma f}(x)$ \eqref{eq:prox_WMAE} and prove its correctness in \cref{theorem:proof_of_alg}.
To this end, we assume that the points $\data_i$ have been sorted and duplicates have been removed and their weights added.
As a result, we can assume $\data_1 < \data_2 < \ldots < \data_N$.
Moreover, we assume $\gamma > 0$, $N \ge 1$ and $w_i > 0$ for all $i = 1, \ldots, N$.
Summands with $w_i = 0$ can obviously be dropped from the sum in \eqref{eq:prox_WMAE}.

We divide the real line into the intervals
\begin{equation}
	I_1
	\coloneqq
	(-\infty,\data_1]
	,
	\quad
	I_i
	\coloneqq 
	\interval{\data_{i-1}}{\data_i}
	\text{ for }
	i = 2, \ldots, N
	\quad
	\text{and}
	\quad
	I_{N+1}
	\coloneqq
	[\data_N,\infty)
	,
\end{equation}
which overlap in the given data points.
It is also useful to set $\data_0 \coloneqq -\infty$ and $\data_{N+1} \coloneqq \infty$.
We further introduce the forward and reverse cumulative weights as
\begin{equation}
	\fwdweight_i
	\coloneqq 
	\sum_{k=1}^i \weight_k
	,
	\quad 
	\revweight_i 
	\coloneqq 
	\sum_{k=i}^N \weight_k
	,
	\quad 
	i = 1, \ldots, N
	.
\end{equation}
We extend these definitions by setting $\fwdweight_0 \coloneqq 0$, $\fwdweight_{N+1}\coloneqq \fwdweight_N$ and $\revweight_{N+1} \coloneqq \revweight_{N+2} \coloneqq 0$. 
We therefore have $\revweight_i = \fwdweight_N - \fwdweight_{i-1}$ for all $i = 1, \ldots, N+2$.
Using this notation, we can rewrite the derivative of $f$ as
\begin{equation}
	\label{eq:derivative_of_mwae}
	f'(y) 
	= 
	\fwdweight_{i-1} - \revweight_i
	\quad 
	\text{for }
	y \in \interior I_i = (\data_{i-1}, \data_i)
	,
	\quad
	i = 1,\ldots,N+1
	.
\end{equation}
This formula reflects the fact that $f$ is piecewise linear and convex since $\fwdweight_{i-1} - \revweight_i$ is monotone increasing with $i$.
Moreover, $f'(y) = \fwdweight_0 - \revweight_1 = - \sum_{k=1}^N \weight_k < 0$ holds for all $y \in \interior I_1$ (points to the left of smallest data point $\data_1$), and $f'(y) = \fwdweight_N - \revweight_{N+1} = \sum_{k=1}^N \weight_k > 0$ holds for all $y \in \interior I_{N+1}$ (points to the right of the largest data point $\data_N$).
At $y = d_i$, $i = 1, \ldots, N$, $f$ is non-differentiable but we can specify its subdifferential, which is
\begin{equation}
	\label{eq:subdifferential_of_mwae}
	\partial f(d_i)
	=
	[\fwdweight_{i-1} - \revweight_i, \fwdweight_i - \revweight_{i+1}]
	,
	\quad
	i = 1, \ldots, N
	.
\end{equation}

The objective
\begin{equation}
	\label{eq:objective_prox_WMAE}
	\Phi(y)
	\coloneqq
	\gamma \, f(y) + \frac{1}{2} (y - x)^2
	=
	\gamma \sum_{i=1}^N \weight_i \, \abs{y - \data_i} + \frac{1}{2} (y - x)^2
\end{equation}
of \eqref{eq:prox_WMAE} is piecewise quadratic and strongly convex.
Its derivative is thus strongly monotone and it satisfies 
\begin{equation}
	\label{eq:limiting_conditions}
	\Phi'(y) < 0 
	\quad
	\text{for all }
	y < \min\{\data_1, x\} 
	\quad
	\text{and}
	\quad
	\Phi'(y) > 0
	\text{ for all }
	y > \max\{\data_N, x\}
	.
\end{equation}
Consequently, the unique minimizer of $\Phi$ lies between these bounds.

\begin{figure}
	\centering
	\begin{tikzpicture}
		\begin{axis}[
			width = 0.6*\linewidth,
			grid = major,
			axis lines = middle,
			axis line style = ->,
			yticklabels = {},
			xticklabels = {},
			xlabel = $y$,
			ylabel = $\partial\obj(y)$,
			xmin=-0.8,xmax=0.55,
			ymin=-1.6,ymax=1.3,
			legend style = {at = {(0.9,0.1)}, anchor=south east,nodes=right},
			]
		\addplot[draw = TolVibrantMagenta,very thick] coordinates {(-0.9,-0.8) (-0.4,-0.3) (-0.4,0.1) (-0.2,0.3) (-0.2,0.4) (0.3,1) (0.3,1.2) (0.6,1.5)};
		\addlegendentry{$\partial\Phi_1(y)$}
		\addplot[draw = TolVibrantTeal,very thick] coordinates {(-0.9,-1.5) (-0.4,-1) (-0.4,-0.8) (-0.2,-0.5) (-0.2, -0.3) (0.3, 0.2) (0.3, 0.35) (0.6,0.65)};
		\addlegendentry{$\partial\Phi_2(y)$}
    \addplot[draw=none,mark=o] coordinates {(-0.4,0) (-0.2,0) (0.3,0)};
    \addplot[draw=none,mark=*,mark options={scale=0.5,}] coordinates {(-0.4,0) (0.1,0)};
    \node[below] at (axis cs: -0.4,0) {$d_1 {{}={}} \color{TolVibrantMagenta}{y_1^*}$};
    \node[below,] at (axis cs: 0.1,0) {$\color{TolVibrantTeal}{y_2^*}$};
    \node[below] at (axis cs: -0.2,0) {$d_2$};
    \node[below] at (axis cs: 0.3,0) {$d_3$};
		\end{axis}
	\end{tikzpicture}
	\caption{Visualization of the subdifferential $\partial\obj$ of the objective $\obj$. In the first (upper) case, $y_1^* = d_1$ holds. In the second (lower) case, we have $y_2^* \in \interior I_3$.}
	\label{figure:SubdifferentialPhi}
\end{figure}
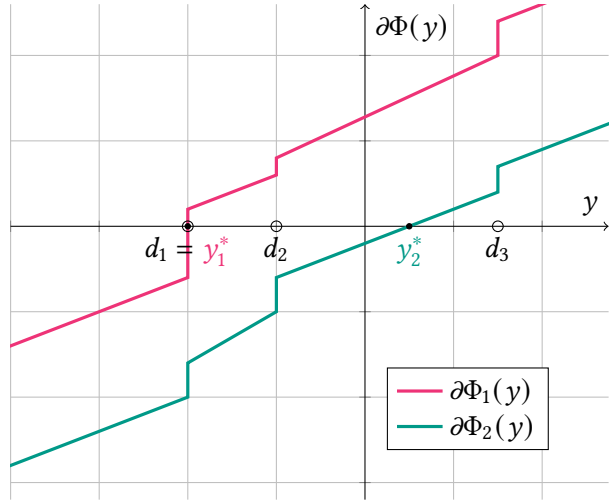

The idea to finding the unique minimizer $y^*$ of \eqref{eq:objective_prox_WMAE} is to locate the smallest index~$1 \le k \le N+1$ such that $y^* \le \data_k$ holds, \ie, the nearest data point to the right of $y^*$.
In other words, we need to find $1 \le k \le N+1$ such that
\begin{subequations}
	\label{eq:determining_the_index}
	\begin{align}
		\label{eq:determining_the_index:1}
		\lim_{y \searrow \data_{k-1}} \Phi'(y)
		&
		=
		\gamma \, (\fwdweight_{k-1} - \revweight_k) + \data_{k-1} - x
		<
		0
		\\
		\label{eq:determining_the_index:2}
		\text{and}
		\quad
		\lim_{y \searrow \data_k} \Phi'(y)
		&
		=
		\gamma \, (\fwdweight_k - \revweight_{k+1}) + \data_k - x
		\ge
		0
	\end{align}
\end{subequations}
holds.
Now we can distinguish two cases: $y^* = d_k$ and $y^* < d_k$.
The first case applies if and only if 
\begin{equation}\label{eq:condition_minimizer_at_datapoint}
	\gamma \, (\fwdweight_{k-1} - \revweight_k) + \data_k - x 
	\le 
	0
	.
\end{equation}
Otherwise, $y^*$ lies in $\interior I_k$ and thus it is the unique minimizer $x - \gamma \, (\fwdweight_{k-1} - \revweight_k)$ of the locally quadratic objective $\obj$.
In either case, once the index~$1 \le k \le N+1$ has been identified, $y^*$ is given by
\begin{equation}
	\label{eq:characteriztion_minimizer_prox_WMAE}
	y^* 
	= 
	\min\{d_k, x - \gamma \, (\fwdweight_{k-1} - \revweight_k) \}
	.
\end{equation}
Both cases are also depicted in \cref{figure:SubdifferentialPhi}.

The considerations above lead to \cref{algorithm:prox_WMAE}.
In our implementation, we evaluate \eqref{eq:prox_WMAE_index_condition} for all $k$ simultaneously and benefit from the quantities being monotone increasing with $k$ when finding the first non-negative entry.
\begin{algorithm}
	\caption{Evaluation of \eqref{eq:prox_WMAE}, the proximal map of the weighted mean absolute error.}
	\label{algorithm:prox_WMAE}
	\begin{algorithmic}[1]
		\Require data points $\data_1 < \data_2 < \ldots < \data_N \in \R$, $N \ge 1$, and $\data_0 \coloneqq -\infty$, $\data_{N+1} \coloneqq \infty$
		\Require weight vector $\weight \in \R^N$ with entries $\weight_i > 0$
		\Require prox parameter $\gamma > 0$ and point of evaluation $x \in \R$
		\Ensure $y = \prox{\gamma f}(x)$, the unique solution of \eqref{eq:prox_WMAE}
		\State Find the smallest index $1 \le k \le N+1$ that satisfies 
		\label{step:determining_the_index}
		\begin{equation}
			\label{eq:prox_WMAE_index_condition}
			\gamma \, (\fwdweight_k - \revweight_{k+1}) + \data_{k} - x
			\ge
			0
		\end{equation}
		\State \Return $y \coloneqq \min \paren[auto]\{\}{d_k, x - \gamma \, (\fwdweight_{k-1} - \revweight_k)}$
		\label{step:determining_the_solution}
	\end{algorithmic}
\end{algorithm}

Let us prove the correctness of \cref{algorithm:prox_WMAE}.
\begin{theorem}
	\label{theorem:proof_of_alg}
	Under the assumptions stated in \cref{algorithm:prox_WMAE}, it returns $y^* = \prox{\gamma f}(x)$, the unique solution of \eqref{eq:prox_WMAE}.
\end{theorem}
\begin{proof}
	Let $1 \le k \le N+1$ be the index found in \cref{step:determining_the_index}.
	First suppose $2 \le k \le N$.
	Then $d_{k-1}$ and $d_k$ are both finite, and \eqref{eq:determining_the_index}, \eqref{eq:prox_WMAE_index_condition} imply
	\begin{equation*}
		\max \partial \Phi(\data_{k-1})
		=
		\lim_{y \searrow \data_{k-1}} \Phi'(y)
		<
		0
		\quad
		\text{and}
		\quad
		\max \partial \Phi(\data_k)
		=
		\lim_{y \searrow \data_k} \Phi'(y)
		\ge
		0
		.
	\end{equation*}
	Owing to the properties of the subdifferential of strongly convex functions, there exists a unique point $y^* \in (d_{k-1},d_k]$ such that $0 \in \partial \Phi(y^*)$, \ie, $y^*$ is the unique minimizer of \eqref{eq:prox_WMAE}.
	This point either belongs to $\interior I_k$, or else $y^* = d_k$ holds.
	In the first case, $\Phi$ is differentiable, so that $\Phi'(y^*) = 0$ holds, yielding
  \begin{equation*}
    x - \gamma \, (\fwdweight_{k-1}-\revweight_k) 
    = 
    y^* 
    < 
		\data_k
    .
  \end{equation*}
	Otherwise, we have $y^* = d_k$, and $0 \in \partial \obj(\data_k)$ implies
	\begin{equation}
		\data_k 
    \le 
		x - \gamma \, (\fwdweight_{k-1} - \revweight_k)
		.
	\end{equation}
	In either case, the unique solution $y^*$ of \eqref{eq:prox_WMAE} is determined by 
	\begin{equation*}
		y^* 
		= 
		\min \paren[auto]\{\}{d_k, x - \gamma \, (\fwdweight_{k-1} - \revweight_k)}
		, 
	\end{equation*}
	which is the quantity returned in \cref{step:determining_the_solution}.

	It remains to verify the marginal cases $k = 1$ and $k = N+1$. 
	In case $k = N+1$, we have $\lim_{y \searrow \data_N} \Phi'(y) < 0$ due to the minimality of $k$. 
	Hence, 
	\begin{equation*}
		\data_N 
		< 
		y^* 
		= 
		x - \gamma \, (\fwdweight_N - \revweight_{N+1}) 
		< 
		\data_{N+1} 
		= 
		\infty
		.
	\end{equation*}
	A similar reasoning applies in the case $k = 1$.
\end{proof}

We provide an efficient and vectorized \python implementation of \cref{algorithm:prox_WMAE} in \cite{BaumgaertnerHerzogSchmidtWeiss:2022:2}. 
It allows the simultaneous evaluation of \eqref{eq:prox_WMAE} for multiple values of $x$, provided that each instance of \eqref{eq:prox_WMAE} has the same number~$N$ of data points.
The weights~$\weight_i$ and data points~$\data_i$ as well as the prox parameter~$\gamma$ may vary between instances.
The discussion so far assumed positive weights for simplicity, but the case $\weight_i = 0$ is a simple extension and it is allowed in our implementation.
This is convenient in order to simultaneously solve problem instances which differ with respect to the number of data points~$N$.
In this case, we can easily pad all instances to the same number of data points using zero weights.
In addition, data points are allowed to be duplicate, \ie, we only require $\data_1 \le \data_2 \le \ldots \le \data_N \in \R$, $N \ge 1$.
Notice that since the data points are assumed to be sorted, finding the index in \cref{step:determining_the_index} is of complexity $\log N$.

\section{Structure of \texorpdfstring{$\prox{\gamma f}$}{prox gamma f}}
\label{section:Structure}

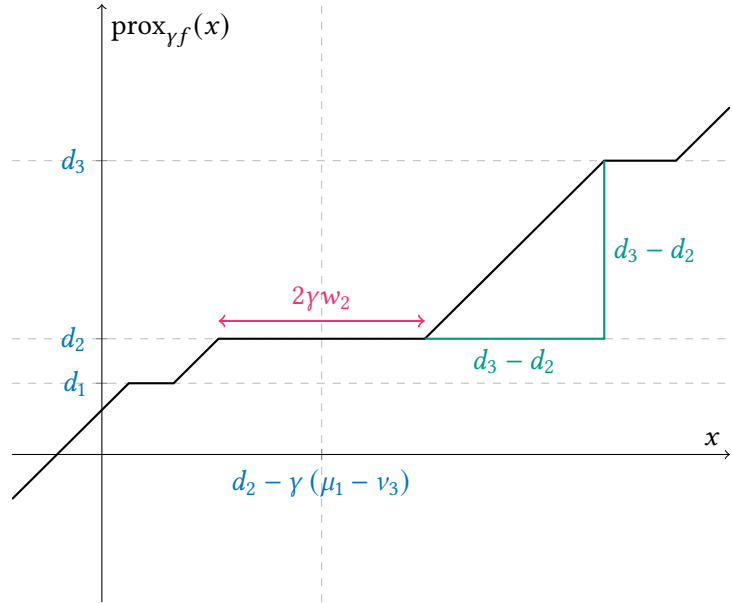
\begin{figure}
	\centering
	\begin{tikzpicture}[outer sep = 0pt]
		\pgfplotstableread{plots/stair_values.csv}{\loadedtable}
		\pgfplotstablegetelem{2}{[index] 0}\of{\loadedtable}
		\begin{axis}[
			width = 0.7*\linewidth,
			grid = major,
			grid style = {dashed},
			axis lines = middle,
			axis line style = ->,
			axis equal = true,
			xlabel = $x$,
			ylabel = $\prox{\gamma f}(x)$,
			xtick = {1.225},
			xticklabels = {$\data_2 - \gamma \, (\fwdweight_1 - \revweight_3)$},
			tick label style = {color = TolVibrantBlue},
			ytick = {0.4,0.65,1.65},
			yticklabels = {$\data_1$,$\data_2$,$\data_3$},
			]
			\addplot[color = black,thick]
				coordinates{(-0.5,-0.25) (0.15,0.4) (0.4,0.4) (0.65,0.65) (1.8,0.65)(2.8,1.65) (3.2,1.65) (3.5,1.95)};
			\draw[<->,TolVibrantMagenta,thick] ({axis cs:0.65,0.75}) -- ({axis cs:1.8,0.75}) node[midway,above]{$2\gamma\weight_2$};
			\draw[TolVibrantTeal,thick] ({axis cs:1.8,0.65}) 
				-- ({axis cs:2.8,0.65}) node[midway,below]{$\data_3-\data_2$}
				-- ({axis cs:2.8,1.65}) node[midway,right]{$\data_3-\data_2$};
		\end{axis}
		\coordinate (a) at (2,2);
		\coordinate (b) at (0,0);
	\end{tikzpicture}
	\caption{Example of the proximal map $\prox{\gamma f}$ in case $N = 3$.}
	\label{figure:stairplot}
\end{figure}

In this section, we briefly discuss the structure of the map $x \mapsto \prox_{\gamma f}(x)$.
Since it generalizes the soft-thresholding operation \eqref{eq:prox_WMAE_1D_solution}, it is not surprising that we obtain a graph which features a staircase pattern. 
An illustrative plot for certain choices of weights $\weight_i$, data points $\data_i$, and prox parameter $\gamma > 0$ is shown in \cref{figure:stairplot}. 
Each of the $N$ distinct data points provides one plateau in the graph.

Two alternating regimes occur for $y^* = \prox{\gamma f}(x)$, as $x$ ranges over $\R$.
First, when $y^* \in \interior I_k$ holds, then \eqref{eq:characteriztion_minimizer_prox_WMAE} implies that $y^*$ is an affine function of $x$ with slope $1$.
This is the case for $x$ whose associated index~$k$ is constant, \ie, 
\begin{equation*}
	x
	\in 
	\gamma \, (\fwdweight_k - \revweight_{k+1})
	+
	[\data_k, \data_{k+1}]
	.
\end{equation*}
As $x$ increases beyond the upper bound of the interval above, $y^*$ enters a constant regime which applies to
\begin{equation*}
	x
	\in
	\gamma \, (\fwdweight_{k-1} - \revweight_{k+1})
	+
	\data_k 
	+ [- \gamma \weight_k, \gamma \weight_k]
	.
\end{equation*}
Notice that the case $N = 1$ reduces to the soft-thresholding map \eqref{eq:prox_WMAE_1D_solution} with only one plateau.

\section{Application To Image Denoising}\label{section:ImageDenoising}

In this section, we present the application of \cref{algorithm:prox_WMAE} to a classical (ROF) total-variation image denoising problem going back to \cite{RudinOsherFatemi:1992:1}. 
Given noisy image data $\imageData_{i,j}$ of size $D_1 \times D_2$, we seek an image $\imageSolution_{i,j}$ of the same dimension which solves
\begin{multline}
	\label{eq:imageDenoising}
	\text{Minimize} 
	\quad
	\ROFobj(\imageSolution)
	\coloneqq
	\frac{1}{2} \sum_{i=1}^{D_1} \sum_{j=1}^{D_2} (\imageSolution_{i,j}-\imageData_{i,j})^2 
	\\
	+ 
	\TVweight \sum_{i=1}^{D_1-1} \sum_{j=1}^{D_2} \abs{\imageSolution_{i+1,j} - \imageSolution_{i,j}}
	+ 
	\TVweight \sum_{i=1}^{D_1} \sum_{j=1}^{D_2-1} \abs{\imageSolution_{i,j+1}-\imageSolution_{i,j}}
	,
	\quad
	u
	\in \R^{D_1 \times D_2}
	.
\end{multline}
Well-known solution approaches to \eqref{eq:imageDenoising} include the primal-dual hybrid gradient method \cite{ChambollePock:2011:1} and the split Bregman iteration \cite{GoldsteinOsher:2009:1}.
The latter requires the solution of a Laplacian problem for $\imageSolution_{i,j}$ in each iteration.
A simpler approach, considered in \cite{LiOsher:2009:1}, is to partition the unknowns in \eqref{eq:imageDenoising} into two disjoint subsets, according to a checkerboard pattern.
In this case, problem \eqref{eq:imageDenoising} with only one subset of unknowns decouples into independent problems, each of which is of type \eqref{eq:prox_WMAE} with weights $\weight_i = 1$ and $\gamma = \beta$ and can be solved efficiently and in parallel using \cref{algorithm:prox_WMAE}.
We give further implementation details below.

Alternating over both subsets of unknowns, one obtains a block-coordinate descent method as proposed in \cite[Sec.~3]{LiOsher:2009:1}.
Unfortunately, such a method does not necessarily converge towards the global minimizer for non-smooth objectives; see for instance \cite[Sec.~2]{FriedmanHastieHoeflingTibshirani:2007:1}.
Therefore, \cite{LiOsher:2009:1} proposed to restart the block-coordinate descent algorithm using a random perturbation of the final iterate upon stagnation.

We depart from this restarting strategy in the following way.
Upon stagnation of the block-coordinate descent method, we evaluate the steepest descent direction by orthogonally projecting (\wrt the Euclidean norm $\norm{\cdot}$) the zero vector onto the subdifferential of the objective $\ROFobj$ from \eqref{eq:imageDenoising}:
\begin{equation}
  \label{eq:subgradient_projection}
  \direction 
	= 
	-\proj{\partial \ROFobj(\imageSolution)}(\bnull) 
	= 
	- \argmin_{s \in \partial \ROFobj(\imageSolution)} \norm{s}^2
	.
\end{equation}
Due to the structure of the subdifferential of the absolute value function $\abs{\cdot}$, this amounts to solving a quadratic optimization problem with $D_1 D_2$ unknowns and sparse linear equality as well as bound constraints, which describe the condition $s \in \partial \ROFobj(\imageSolution)$. 
We employ the QP solver \osqp (\cite{StellatoBanjacGoulartBemporadBoyd:2020:1}, \url{https://github.com/osqp/osqp}) for the purpose of solving \eqref{eq:subgradient_projection}.

The steepest descent direction \eqref{eq:subgradient_projection} is used in the following way in \cref{algorithm:checkerboard_denoising}.
First, it serves as a perturbation direction upon stagnation, in contrast to the random perturbation proposed in \cite{LiOsher:2009:1}; see \cref{step:step_size}.
Second, the norm $\norm{\direction}$ can serve as a stopping criterion; see \cref{step:stopping}.

\begin{algorithm}
	\caption{Checkerboard scheme to approximately solve \eqref{eq:imageDenoising}.}
	\label{algorithm:checkerboard_denoising}
	\begin{algorithmic}[1]
		\Require tolerances $\tolinner$, $\tolouter$
		\Ensure $\imageSolution^k$, an approximate solution of \eqref{eq:imageDenoising}
		\State $k \coloneqq 0$
		\Repeat 
		\Repeat 
		\State $\imageSolution_w^{k+1} \coloneqq \argmin\limits_{\imageSolution_w} \ROFobj \paren[big](){\imageSolution_b^k + \imageSolution_w} $ using \cref{algorithm:prox_WMAE}
		\label{step:white_step}
		\State $\imageSolution_b^{k+1} \coloneqq \argmin\limits_{\imageSolution_b} \ROFobj \paren[big](){\imageSolution_b + \imageSolution_w^{k+1}}$ using \cref{algorithm:prox_WMAE}
		\label{step:black_step}
		\State Set $k \coloneqq k+1$
		\label{step:increase_iter}
		\Until{$\norm{\imageSolution^k - \imageSolution^{k-1}} \le \tolinner$}
		\label{step:stopping_inner}
		\State $\direction \coloneqq -\proj{\partial \ROFobj(\imageSolution^k)}(\bnull)$
    \label{step:projection}
		\If{$\norm{\direction} > \tolouter$}
		\State Choose $\alpha > 0$ with $\ROFobj(\imageSolution^k + \alpha \direction) < \ROFobj(\imageSolution^k)$
		 \label{step:step_size}
		\State Set $\imageSolution^{k+1} \coloneqq \imageSolution^k + \alpha \direction$
		\State Set $k \coloneqq k+1$
		\EndIf
		\label{step:end}
		\Until{$\norm{\direction} \le \tolouter$}
		\label{step:stopping}
	\end{algorithmic}
\end{algorithm}

In \cref{step:black_step,step:white_step}, we use subvector indexing.
That is, $\imageSolution_w$ refers to the subvector of $u$ with \enquote{white} indices copied and \enquote{black} indices zeroed.
The roles are reversed for $\imageSolution_b$.
Consequently, $u = \imageSolution_b + \imageSolution_w$ holds.
Subvector indexing can be conveniently done in \python using logical indexing.

Notice that \cref{step:black_step,step:white_step} require the evaluation of a function of the \eqref{eq:prox_WMAE} for many arguments in parallel, where we benefit from our vectorized implementation of \cref{algorithm:prox_WMAE}.
All norms in \cref{theorem:proof_of_alg} are Frobenius norms for matrices.
\cref{step:step_size} is implemented using a backtracking strategy starting with initial step size $\alpha = 0.5$, which is then halved until the condition $\ROFobj(\imageSolution^k + \alpha \direction) < \ROFobj(\imageSolution^k)$ is met.

We wish to emphasize that we do not propose \cref{algorithm:checkerboard_denoising} as a novel solver for image denoising problems.
We rather consider it here as a source of problems of type \eqref{eq:prox_WMAE}, which can be solved efficiently by our proposed \cref{algorithm:prox_WMAE}.
In practice, \cref{algorithm:checkerboard_denoising} can also be used effectively as a preliminary solver stage for \eqref{eq:imageDenoising}, before one switches to, \eg, the split Bregman or Chambolle-Pock iteration.

\begin{figure}[hptb]
	\centering
	\begin{subfigure}[t]{0.3\textwidth}
    \centering
    \includegraphics[width=0.9\linewidth]{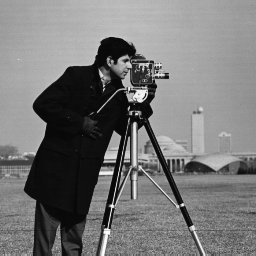}
	  \caption{Noisy cameraman image.}
    \label{figure:cameraman_noise}
	\end{subfigure}%
	\hfill
	\begin{subfigure}[t]{0.3\textwidth}
    \centering
    \includegraphics[width=0.9\linewidth]{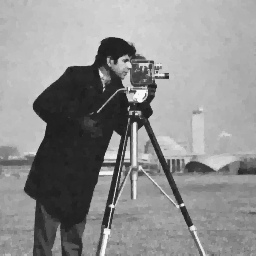}
		\caption{Approximate solution obtained with \cref{algorithm:checkerboard_denoising}.}
    \label{figure:cameraman_checkerboard} 
	\end{subfigure}%
	\hfill
	\begin{subfigure}[t]{0.3\textwidth}
    \centering
    \includegraphics[width=0.9\linewidth]{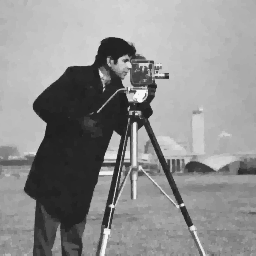}
    \caption{Exact solution of \eqref{eq:imageDenoising} obtained by a split Bregman iteration.}
    \label{figure:cameraman_result}
	\end{subfigure}
	\caption{Numerical results to obtain timings for \cref{algorithm:prox_WMAE}.}
  \label{figure:cameraman_comparison_cb_exact}
\end{figure}

For verification purposes, we show the outcome of \cref{algorithm:checkerboard_denoising}.
As a test case, we choose the well-known cameraman image of size $D_1 = D_2 = 256$. 
We add zero-mean Gaussian noise with standard deviation~$\sigma = 50$ independently pixel by pixel.
We then truncate the values to the range $[0,255]$ to obtain the noisy image shown in \cref{figure:cameraman_noise}. 
We apply \cref{algorithm:checkerboard_denoising} to the image denoising problem \eqref{eq:imageDenoising} with parameter $\TVweight = 10$.
The inner tolerance is set to $\tolinner = 10^{-4}$.
We observe that even for a coarse outer tolerance of $\tolouter = 300$, a good reconstruction is obtained; see \cref{figure:cameraman_checkerboard}.
This tolerance is reached after $k = 42$ iterations, of which $37$ are iterations of the loop in \cref{step:white_step}--\cref{step:increase_iter} and $5$ are executions of \cref{step:projection}--\cref{step:end}.
In particular, the subdifferential projection step in \cref{step:projection}, which amounts to solving a quadratic optimization problem, is carried out $5$~times.
For comparison, we include an \enquote{exact} solution of \eqref{eq:imageDenoising} obtained by a split Bregman iteration with very tight tolerances in \cref{figure:cameraman_comparison_cb_exact}.

Each call to \cref{algorithm:prox_WMAE} (\cref{step:white_step,step:black_step} in \cref{algorithm:checkerboard_denoising}) evaluates the proximity operator \eqref{eq:prox_WMAE} in parallel for half the number of pixels, \ie, $2^{15} = \num{32768}$ instances.
Using padding with zero weights for instances of \eqref{eq:prox_WMAE} pertaining to points on the boundary, all instances have $N = 4$ data points $\data_i$ (the current values for all neighbors north, south, east and west) and weights $\weight_i \in \{0,1\}$.
Timing results are reported in \cref{table:timings}.
They were obtained on a laptop with an 8-core Intel~Core~i5 CPU with 1.6~GHz and 16~GiB RAM, running Ubuntu~22.04 and \python~3.10.

\begin{table}
	\centering
	\begin{tabular}{lrrr}
		\toprule
		subroutine 
		& 
		number of calls
		&
		total time
		&
		time per call
		\\
		\midrule
		\cref{algorithm:prox_WMAE} 
		\\
		(\num{32768} solves of \eqref{eq:prox_WMAE} in parallel)
		&
		74
		&
		\SI{10.7}{\second}
		&
		\SI{0.14}{\second}
		\\
		\midrule
		solution of QP (\cref{step:projection}) using \osqp
		&
		5
		&
		\SI{11.7}{\second}
		&
		\SI{2.34}{\second}
		\\
		\bottomrule
	\end{tabular}
	\caption{Timing results.}
	\label{table:timings}
\end{table}

\section{Application to a Deflection Energy Minimization Problem}\label{section:Application_to_Deflection_Energy_Minimization}

In this section, we consider the minimization of the non-smooth energy~$\energy$ of a deflected membrane.
We apply an ADMM (alternating direction method of multipliers) scheme which requires the parallel evaluation of the proximal map \cref{algorithm:prox_WMAE} in each inner loop of the scheme. 

Consider a bounded domain $\Omega \subset \R^2$ occupied by a thin membrane.
In our model, the energy of this membrane is given in terms of its unknown deflection (displacement) function $z \colon \Omega \to \R$ and it takes the following form: 
\begin{equation}\label{eq:energy}
  \energy (z) 
  = 
  \frac{1}{2} \int_\Omega c \, \abs{\nabla z}^2 \d x
  -
  \int_\Omega f z \d x
  + 
	\sum_{i=1}^L \weight_i \int_\Omega \max \paren\{\}{z - \data_i, \, 0} \d x 
  + 
  \frac{1}{2} \int_\Gamma \alpha \, z^2 \d s
  .
\end{equation}
The first term describes the potential energy of the membrane, where $c \in \R$ is the stiffness constant.
The second term accounts for an external area force density~$f$, which we assume to be a non-negative function on $\Omega$.
Consequently, the deflection will be non-negative as well.
The specialty of our model is the third term, which can be interpreted as follows.
The non-negative constants $\data_1, \ldots, \data_L$ serve as thresholds.\footnote{We could also allow $\data_1, \ldots, \data_L$ to be non-negative functions on $\Omega$, with minor modifications in what follows.}
Once the deflection~$z$ at a point in $\Omega$ exceeds any of these thresholds~$\data_i$, an additional downward force of size $\weight_i > 0$ times the excess deflection activates.
Finally, the fourth term models an additional potential energy due to springs along the boundary~$\Gamma$ of the domain with stiffness constant~$\alpha$.
A related problem has been studied in \cite{BostanHanReddy:2005:1} with an emphasis on a~posteriori error analysis and adaptive solution.

Notice that the minimization of the convex, non-smooth energy~$\energy$ among all displacements~$z \in H^1(\Omega)$ corresponds to the weak form of a partial differential equation (PDE), or rather a variational inequality.
In fact, the necessary and sufficient optimality conditions for the minimization of \eqref{eq:energy} amount to
\makeatletter
\begin{multline}\label{eq:PDE_weak_form}
  c
  \int_\Omega \nabla z \cdot \nabla \delta z \d x 
  - 
  \int_\Omega f \, \delta z \d x 
  + 
	\sum_{i=1}^L \weight_i \int_\Omega \chi_{\setDef{x \in \Omega}{z(x) > \data_i}}
  \, \delta z \d x
	\\
  + 
	\sum_{i=1}^L \weight_i \int_\Omega \chi_{\setDef{x \in \Omega}{z(x) = \data_i}}
	\, \max\{\delta z, \, 0\} \d x
  + 
  \int_\Gamma \alpha \, z \, \delta z \d s 
	\ge
	0
	\ltx@ifclassloaded{svjour3}{\\}{\quad}
	\text{for all }
  \delta z \in H^1(\Omega)
  .
\end{multline}
\makeatother
Here, $\chi_A\colon \Omega \to \R$ is the characteristic function of the set~$A$ with values in $\{0,1\}$.  
Provided that the set $\setDef{x \in \Omega}{z(x) = \data_i}$ is of Lebesgue measure zero, the variational inequality \eqref{eq:PDE_weak_form} becomes an equation, whose strong form---using integration by parts---can be seen to be
\begin{equation*}
	\begin{alignedat}{2}
		\label{eq:PDE_strong_form}
		- c \, \laplace z
		+ 
		\sum_{i=1}^L \weight_i \chi_{\setDef{x \in \Omega}{z(x) > \data_i}}
		&
		= 
		f
		&
		&
		\quad
		\text{in }
		\Omega
		,
		\\
		c 
		\frac{\partial z}{\partial n}
		&
		= 
		\alpha z
		&
		&
		\quad 
		\text{ on } 
		\Gamma
		.
	\end{alignedat}
\end{equation*}
From here we also learn that the boundary term in the energy \eqref{eq:energy} leads to boundary conditions of Robin type.

We employ a standard Galerkin approach to numerically discretize the energy \eqref{eq:energy}.
To this end, we replace the displacement space $H^1(\Omega)$ by some finite dimensional subspace of piecewise linear, globally continuous finite element functions defined over a triangulation of $\Omega$.
Denoting the associated nodal basis by $\{\varphi_i\}$, we define the stiffness matrix 
\begin{equation*}
	K_{ij} 
	= 
	\int_\Omega c \nabla \varphi_i \cdot \nabla \varphi_j \d x
	+
	\int_\Gamma \alpha \, \varphi_i \, \varphi_j \d s
\end{equation*}
as well as the lumped mass matrix
\begin{equation*}
	M_{ii}
	=
	\int_\Omega \varphi_i \d x
\end{equation*}
with $M_{ij} = 0$ when $i \neq j$.
The reason for using mass lumping to approximate all area integrals in \eqref{eq:energy} which do not involve derivatives is to translate the pointwise maximum operator into a coefficientwise one.
This technique is crucial for an efficient numerical realization and has been used before, \eg, in \cite{WachsmuthWachsmuth:2011:3,CasasHerzogWachsmuth:2012:3}.

We continue to use $z$ to denote the nodal values of the finite element discretization of the deflection and thus obtain the discrete energy
\begin{equation}\label{eq:energy_discrete_max}
	\energyh(z)
	= 
	\frac{1}{2} z^\transp K z 
	- 
	f^\transp M z
	+
	\sum_{i=1}^L \weight_i \bone^\transp M \, \max\{z - \data_i, \, 0\}
	.
\end{equation}
Here, $\bone$ denotes the vector of all ones and the $\max$ operation is understood coefficientwise. 

The minimization of the convex but non-smooth energy \eqref{eq:energy_discrete_max} is not straightforward.
Our method of choice here is an ADMM scheme, which moves the non-smooth terms into a separate subproblem, which can then be efficiently solved using \cref{algorithm:prox_WMAE}.
We refer the reader to \cite{BoydParikhChuPeleatoEckstein:2010:1} and the references therein for an overview on ADMM.

In order to obtain a subproblem of type \eqref{eq:prox_WMAE}, we use the identity $2 \max\{a, b\} = a + b + \abs{a-b}$ and arrive at
\begin{equation}\label{eq:energy_discrete}
	\energyh(z)
	= 
	\frac{1}{2} z^\transp K z 
	- 
	\tilde f^\transp M z
	+
  \frac{1}{2} \sum_{i=1}^L \weight_i \bone^\transp M \, \abs{z - \data_i}
  + 
	C
\end{equation}
where $C \coloneqq \sum_{i=1}^L \weight_i \data_i \bone^\transp M \bone$ is a constant and $\tilde f \coloneqq f - \frac{1}{2} \sum_{i=1}^L \weight_i \bone$ is a modification of the force vector.
Moreover, the absolute value $\abs{\,\cdot\,}$ is understood coefficientwise.
Dropping the constant~$C$, introducing a second variable~$y$, the constraint $z = y$ and associated (scaled) Lagrange multiplier~$\mu$ gives rise to the augmented Lagrangian
\begin{multline}\label{eq:AugmentedLagrangian}
	\cL_\rho(z,y,\mu) 
	\\
	= 
	\frac{1}{2} z^\transp K z 
	- 
	\tilde f^\transp M z
	+
  \frac{1}{2} \sum_{i=1}^L \weight_i \bone^\transp M \, \abs{y - \data_i \bone}
	+ 
	\frac{\rho}{2}
	\norm{z - y + \mu}_M^2
	,
\end{multline}
where $\norm{\,\cdot\,}_M^2$ is the squared norm induced by the positive diagonal matrix~$M$.

An ADMM computes iterates $\sequence{z}{k}$, $\sequence{y}{k}$, $\sequence{\mu}{k}$ according to the following update scheme 
\begin{subequations}
	\label{eq:deflection_ADMM}
	\begin{align}
		z^{(k+1)}
		&
		\coloneqq
		\argmin_z \cL_\rho(z, \sequence{y}{k}, \sequence{\mu}{k})
		,
		\label{eq:deflection_z_update}
		\\
		\sequence{y}{k+1}
		&
		\coloneqq
		\argmin_y
		\cL_\rho(z^{(k+1)}, y, \sequence{\mu}{k})
		\label{eq:deflection_y_update}
		,
		\\
		\mu^{(k+1)} 
		&
		\coloneqq
		\sequence{\mu}{k} + z^{(k+1)} - y^{(k+1)}
		\label{eq:deflection_mu_update}
		.
	\end{align}
\end{subequations}
For the problem at hand, \eqref{eq:deflection_z_update} amounts to the solution of the discretized Poisson-like problem
\begin{equation*}
	K z + \rho M z
	=
	M \, \paren[big](){\tilde f + \rho \, (\sequence{y}{k} - \sequence{\mu}{k})}
	.
\end{equation*}
The $y$-update \eqref{eq:deflection_y_update} on the other hand can be cast as the problem
\begin{equation}\label{eq:prox_in_deflection}
	\text{Minimize}
	\quad
	\frac{1}{2} \sum_{i=1}^L \weight_i \bone^\transp M \, \abs{y - \data_i \bone}
	+ 
	\frac{\rho}{2}
	\norm{\sequence{z}{k+1} - y + \sequence{\mu}{k}}_M^2
	\quad
	\text{\wrt\ }
	y
	.
\end{equation}
Owing to the diagonal structure of the mass matrix, this problem fully decouples, the mass matrix cancels, and we obtain
\begin{equation}\label{eq:prox_in_deflection_componentwise}
	\text{Minimize}
	\quad
	\frac{1}{2} \sum_{i=1}^L 
	\weight_i \, \abs{y_j - \data_i}
	+ 
	\frac{\rho}{2}
	\paren[big](){\sequence{z_i}{k+1} - y_j + \sequence{\mu_i}{k}}^2
	\quad
	\text{\wrt\ }
	y_j
\end{equation}
for each component $y_j$ of $y$.
This problem fits into the pattern of \eqref{eq:prox_WMAE}, and thus can be solved efficiently and simultaneously for all components~$y_j$ by \cref{algorithm:prox_WMAE}. 

\begin{figure}[htbp]
	\centering
	\begin{subfigure}[t]{0.48\textwidth}
		\centering
		\includegraphics[width=0.9\linewidth]{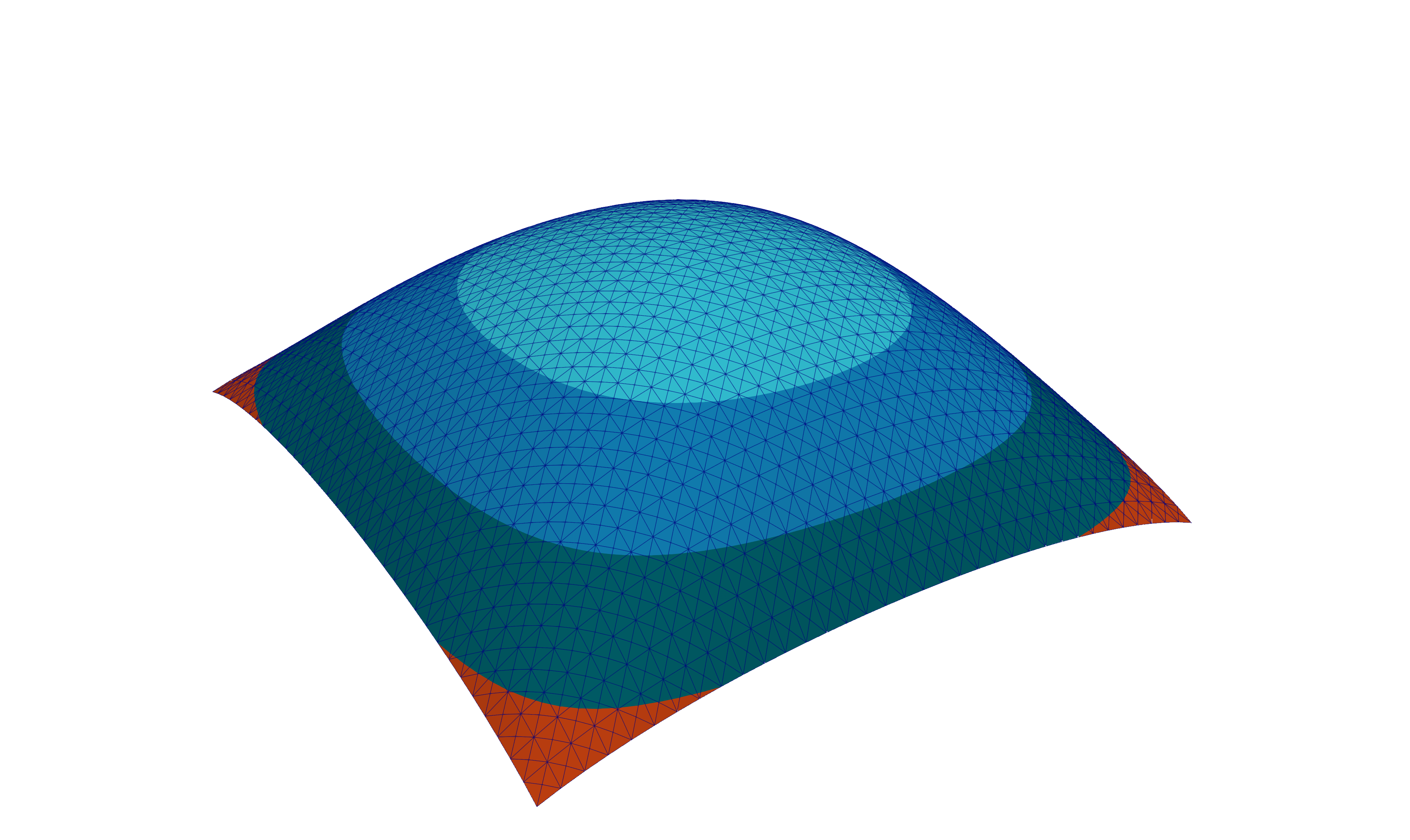}
	\end{subfigure}%
	\hfill
	\begin{subfigure}[t]{0.48\textwidth}
		\centering
		\includegraphics[width=0.9\linewidth]{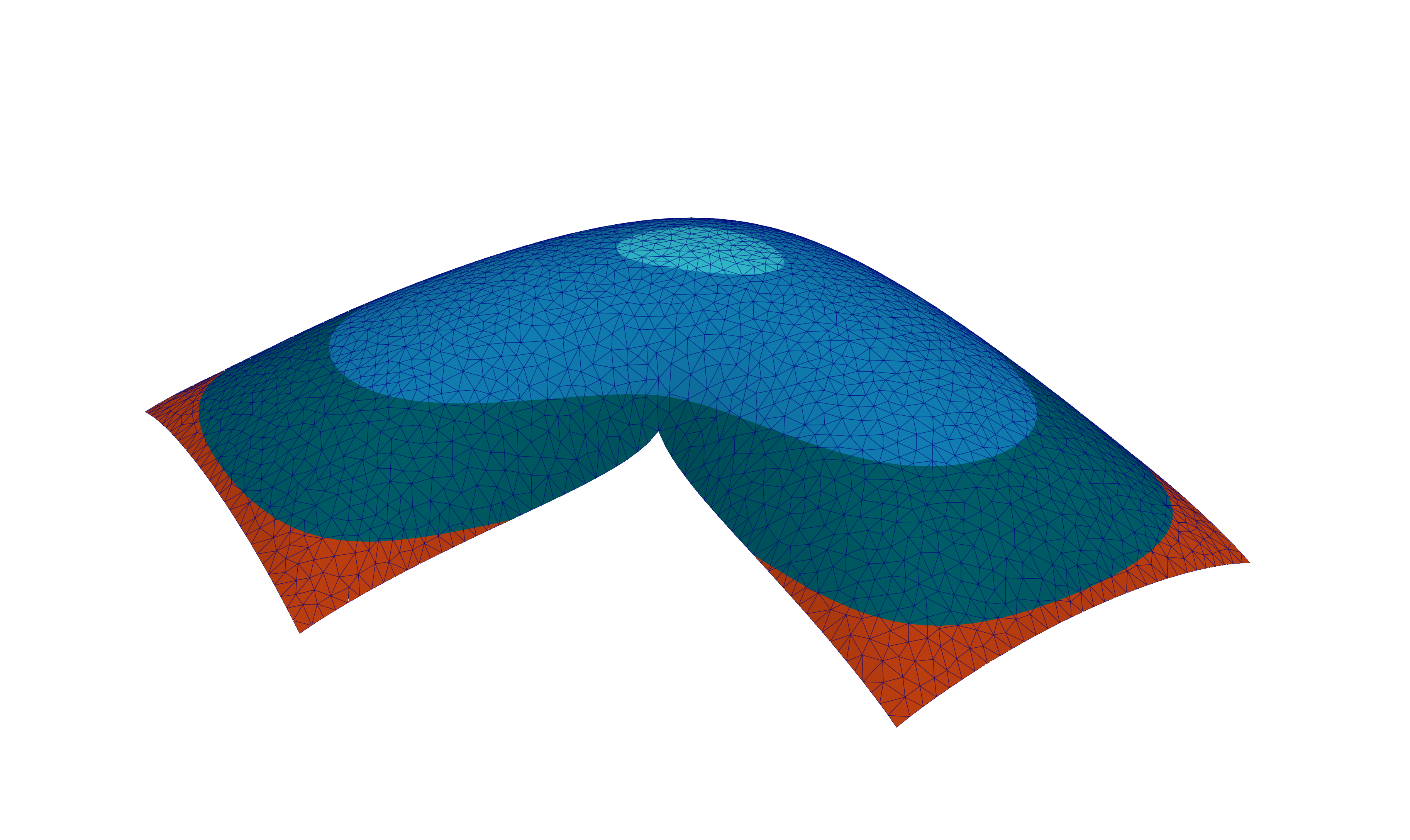}
	\end{subfigure}%
	\hfill
	\caption{Resulting deflections $z$ which minimize the discrete energy functional $\energyh$ defined in \eqref{eq:energy_discrete_max} over two different domains, with data given in \eqref{eq:ExampleConfigurationDeflection}. Deflections were computed with the ADMM scheme \eqref{eq:deflection_ADMM} and stopping criterion \eqref{eq:deflection_stopping_criterion}. The subdomains are colored according to the number of active terms ($z \ge d_i)$ in \eqref{eq:energy_discrete_max}.}
	\label{figure:DeflectionDomains}
\end{figure}

In order to validate the method, we consider two simple, convex domains 
\begin{equation*}
  \Omega_1 \coloneqq (0,1)^2
  ,
	\quad
  \Omega_2
  \coloneqq
  (0,1.1)^2 \setminus (0.6,1.1)^2
  .
\end{equation*}
We use meshes with $2521$~vertices for $\Omega_1$ and $2637$~vertices for $\Omega_2$.
We choose the stiffness constant~$c$, force density~$f$, boundary stiffness constant~$\alpha$, additional downward forces~$\weight_i$ and threshold deflections $\data_i$ to be constant over the domains in each case, as
\begin{equation}\label{eq:ExampleConfigurationDeflection}
  c \equiv 1
  ,
	\quad
  f \equiv 0.5
  ,
	\quad
  \alpha \equiv 10
  ,
	\quad
  \data_i \equiv 0.01 \cdot i
  ,
	\quad
  \weight_i \equiv 0.02
  ,
	\quad
  i = 1, \ldots, 4
  .
\end{equation}
We apply the ADMM scheme \eqref{eq:deflection_ADMM}. 
The iterations are terminated as soon as the stopping criterion
\begin{equation}\label{eq:deflection_stopping_criterion}
	\max \paren[big]\{\}{%
		\norm{\sequence{z}{k+1} - \sequence{z}{k}}_M
		,
		\;
		\norm{\sequence{y}{k+1} - \sequence{y}{k}}_M
		,
		\;
		\norm{\sequence{\mu}{k+1} - \sequence{\mu}{k}}_M
	}
  \le 
  \toldeflection
  \coloneqq 
  10^{-20}
\end{equation}
is satisfied, where we recall that $\norm{\,\cdot\,}_M$ denotes the norm induced by the lumped mass matrix~$M$.
With this setup and $\rho = 100$, the ADMM algorithm required $186$~iterations for the problem on $\Omega_1$ and $278$~iterations in case of $\Omega_2$.
The resulting solutions~$z$ are shown in \cref{figure:DeflectionDomains}. 

